\documentclass[11pt]{article}
\usepackage{amsmath}
\usepackage{amssymb}
\usepackage{amsthm}
\usepackage[dvips]{graphicx,color}

\newtheorem {theorem}{Theorem}[section]

\newtheorem{lemma}[theorem]{Lemma}
\newtheorem{proposition}[theorem]{Proposition}
\newtheorem{corollary}[theorem]{Corollary}

\newtheorem{definition}[theorem]{Definition}
\newtheorem{remark}[theorem]{Remark}

\newtheorem{problem}[theorem]{Problem}
\newtheorem{claim}[theorem]{Claim}

\newcommand {\Heads}[1]   {\smallskip\pagebreak[1]\noindent{\bf
#1{\hskip 0.2cm} }}
\newcommand {\Head}[1]    {\Heads{#1:}}
\newcommand {\proofs}     {\proof}
\newcommand {\prooft}[1]  {\Head{Proof #1}\nopagebreak[2]}
\newcommand  {\QED}    {\def\qedsymbol{$\blacksquare$}\qed}

\def\Q{{\Bbb Q}}
\def\R{{\Bbb R}}

\def\N{{\Bbb N}}
\def\Z{{\Bbb Z}}
\def\P{{\bf{P}}}

\def\eps{\epsilon}
\def\D{{\mathcal D}}
\def\FF{R}

\begin{document}
\title{New coins from old: computing with unknown bias}

\author{Elchanan Mossel \thanks{Supported by a Miller Fellowship}\\
U.C. Berkeley \\ mossel@stat.berkeley.edu \and
Yuval Peres \thanks{Supported in part by NSF Grant DMS-0104073 and by a Miller
  Professorship}
\\ U.C. Berkeley \\
peres@stat.berkeley.edu \\ 
\and With an appendix by Christopher Hillar 
\thanks{This work is supported under a National Science Foundation
  Graduate Research Fellowship.} 
\\ U.C. Berkeley 
\\ chillar@math.berkeley.edu}
\date{\today}
\maketitle

\begin{abstract}
Suppose that we are given a function  
$f : (0,1)  \to (0,1)$  and, for some unknown
 $p \in (0,1)$, a sequence of independent tosses of a $p$-coin
(i.e., a coin with probability $p$ of ``heads'').
For which functions $f$
is it possible to simulate an $f(p)$-coin? $\;$
This question was raised by S. Asmussen and J. Propp.
A simple simulation scheme for the constant function $f(p) \equiv 1/2$
was described by von Neumann (1951);
this scheme can be easily implemented using a finite automaton.
We prove that in general, an $f(p)$-coin can be simulated by a 
finite automaton for all $p \in (0,1)$,
if and only if $f$ is a rational function over $\Q$. 
We also show that if an $f(p)$-coin can be simulated by a pushdown automaton, 
then $f$ is an algebraic function over $\Q$; 
however, pushdown automata can simulate 
$f(p)$-coins for certain non-rational functions such as $f(p)=\sqrt{p}$.
These results complement the work of Keane and O'Brien (1994), who
determined the functions $f$ for which an $f(p)$-coin can be simulated
when there are no computational restrictions on the simulation scheme.


\end{abstract}

\section{Introduction}
Fifty years ago, von Neumann \cite{vN} suggested a method to 
generate unbiased random bits from a sequence of i.i.d.\ biased bits. 
This method can be easily implemented using a finite automaton.


In this paper we study the following generalization.
Let $\D \subset (0,1)$. Suppose that we are 
given a function  $f : \D \to (0,1)$  and, for some unknown
 $p \in \D$, a sequence of  independent tosses of a $p$-coin
(i.e., $\{0,1\}$ valued random variables with mean $p$).
For which functions $f$ is it then possible to simulate 
an $f(p)$-coin? 

The allowed simulation schemes  apply
a stopping rule to independent tosses of a $p$-coin, and then
determine a $\{0,1\}$-valued variable with
mean $f(p)$ as a function of the stopped sequence.
 We emphasize that the scheme cannot depend on $p$.
We are especially interested in simulation schemes that can be implemented
by an automaton that receives the $p$-coin tosses as inputs,
and outputs an $f(p)$-coin; see \S 1.2 for more formal definitions.

A special case of this question was raised in 1991 by S. Asmussen 
(see  \cite{KO}). We learned of the general problem from
J. Propp (personal communication) who emphasized its computational
aspects. The problem was considered in the context
of Markov chain simulation by Glynn and Henderson~\cite{GH}.

In our main result we prove 
\begin{theorem} \label{thm:newauto}
Let $\D \subset (0,1)$ and  $f : \D \to (0,1)$.
Then an $f(p)$-coin for $p \in \D$ can be 
simulated  using a finite automaton from independent tosses of a $p$-coin,
if and only if $f$ is the restriction to $\D$ 
of a rational function  $F$ over $\Q$, such that $0<F(x)<1$ for
$0<x<1$.
\end{theorem}

Later using the result of the appendix we prove:
\begin{theorem} \label{thm:newpush}
Let $f : (0,1) \to (0,1)$. If an $f(p)$-coin can be simulated 
from tosses of a $p$-coin by a pushdown automaton for all $p \in (0,1)$, 
then $f$ is an algebraic function over $\Q$; 
i.e., there exists a non-zero polynomial $P \in \Q[X,Y]$ such that 
$P(f(p),p) = 0$ for all $p \in [0,1]$.
\end{theorem}
In \S 3.2 we describe pushdown automata that can simulate 
$f(p)$-coins for certain non-rational functions such as $f(p)=\sqrt{p}$.

Our results complement the work of
Keane and O'Brien (1994), who considered the simulation
problem without any computational restrictions. They showed that
the  functions $f:\D \to (0,1)$ for which an $f(p)$-coin can be
simulated (in principle) are precisely the constants,
and all continuous functions on $\D$ that satisfy
$\min(f(x),1-f(x)) \ge \min(x,1-x)^n$  for some $n \ge 1$
and all $x \in \D$. 

\subsection{Examples}
Here are some examples - see Figure 1.
\begin{figure}[h] \label{figure:auto1}
\begin{center}
\includegraphics[height=4cm]{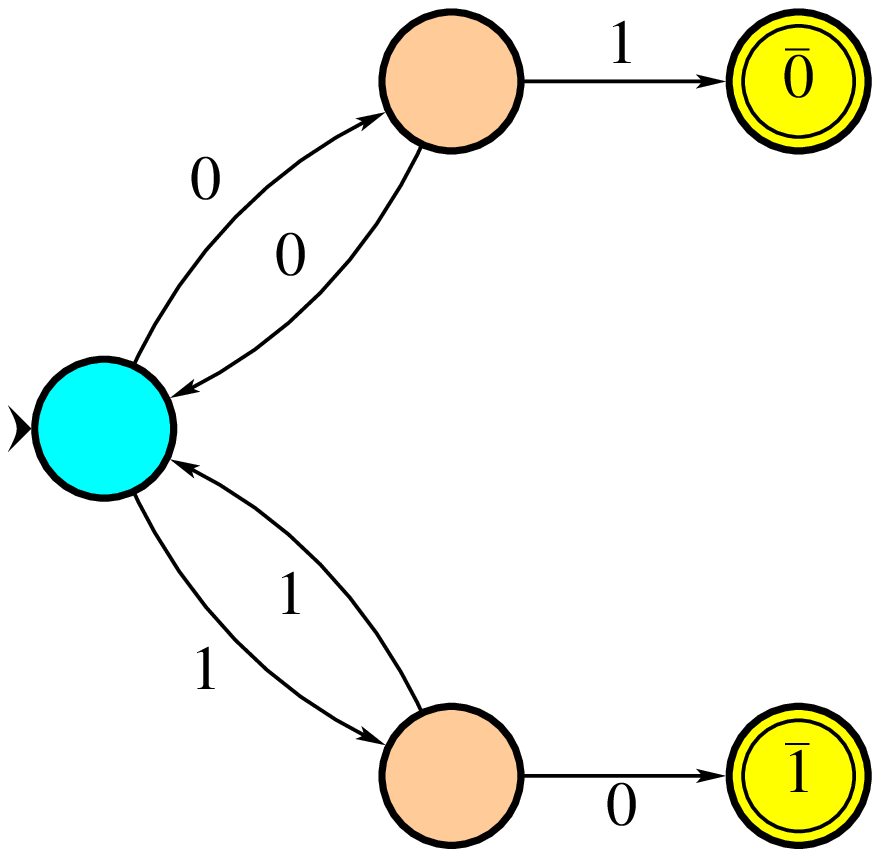} 
\includegraphics[height=4cm]{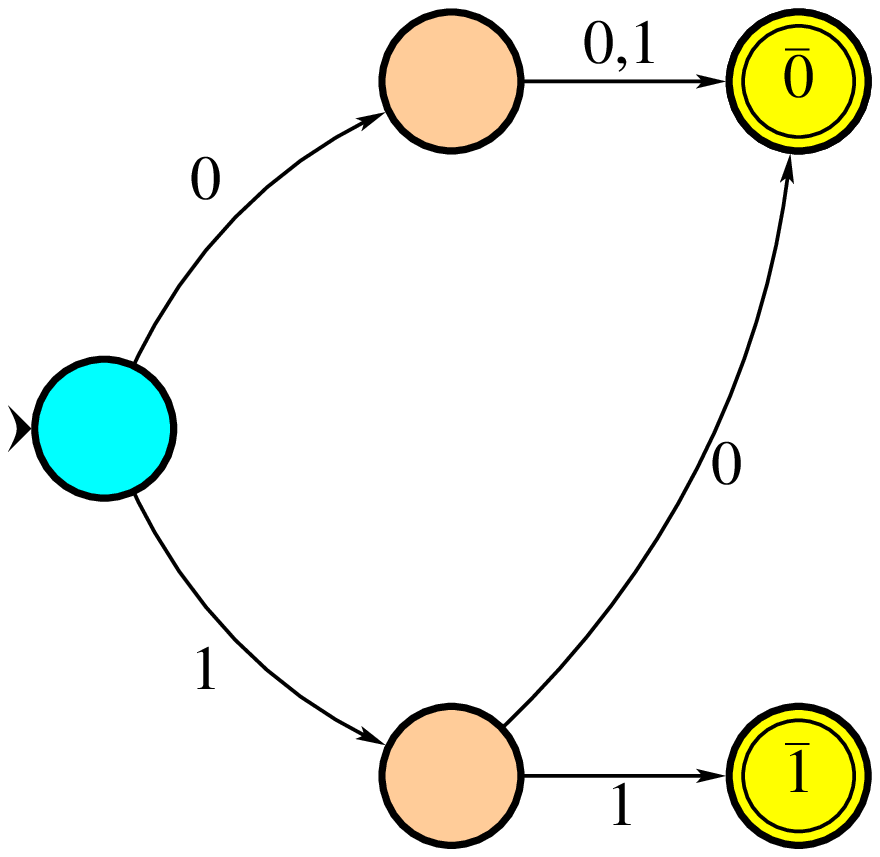}
\includegraphics[height=4cm]{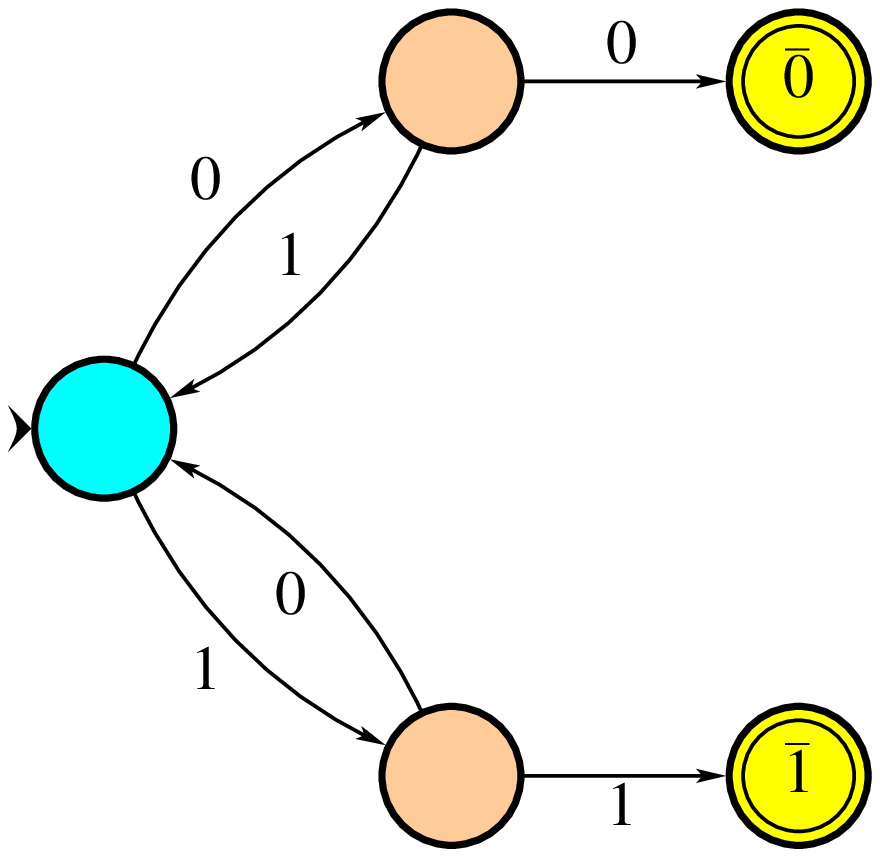}  
\caption{Simulating an unbiased coin, $p^2$ and $\frac{p^2}{p^2 + (1-p)^2}$}
\end{center}
\end{figure}
\begin{itemize}
\item
$f(p) = 1/2$. This is achievable by Von Neumann's trick: toss the $p$ coin twice and let $x$ and
$y$ be the outcome. If $xy = 01$ output $0$, if $xy = 10$ declare $1$;
otherwise, do not declare and toss again. This is the leftmost automaton.
\item
$f(p) = p^2$. Toss the $p$-coin twice. If it's $11$ declare $1$,
otherwise declare $0$. This is the automaton in the middle.
\item
$f(p) = \frac{p^2}{p^2 + (1-p)^2}$. Toss the $p$-coin twice until you
get $00$ or $11$. 
In the first case, declare $0$, in the second declare $1$.
This is the automaton on the right.
\item
$f(p) = \sqrt{p}$. Theorem \ref{thm:auto} implies that  
this function cannot be simulated by a finite automaton.
In Section \ref{sec:push} we construct a pushdown automaton which
simulates $f$.

\item Our main result implies that
there is no infinite set $\D$
such that a $2p$-coin can be simulated by a finite automaton
for all $p \in \D$ from tosses of a $p$-coin.

\end{itemize}

\subsection{Formal definitions}
Denote by $\{0,1\}^{\ast}$ the set of all finite binary strings,
and call any subset of $\{0,1\}^{\ast}$ a {\em language}.
 Say that a language $L$  has the {\em prefix property\/}  if there is
no pair of distinct strings $u,v \in L$ such that 
$u$ is a prefix of $v$.

For any binary string $w$, write
$\P_p[w] = p^{n_1(w)} (1 - p)^{n_0(w)}$, 
where $n_i(w)$ is the number of $i$'s in $w$
and let $\P_p(L)=\sum_{w \in L} \P_p(w)$ for any language $L$.

\begin{definition}
Let $f : \D \to [0,1]$ where $\D \subset (0,1)$.
\begin{itemize}
\item
A {\em simulation of $f$} is a pair of disjoint 
languages $L_0, L_1 \subset \{0,1\}^{\ast}$
such that $L_0 \cup L_1$ has the prefix property 
and for all $p \in \D$ we have $\P_p[L_1] = f(p)$ and
$\P_p[L_0 \cup L_1] = 1$ and $\P_p[L_1] = f(p)$. 
\item
A {\em simulation of $f$ by a finite
  automaton}, is a
simulation $(L_0,L_1)$ of $f$ such that 
there exists a finite automaton
which outputs $0$ for the language $L_0$ and outputs $1$ 
for the language $L_1$. An analogous definition applied to simulation
by a pushdown  automaton or a Turing machine.
\end{itemize}
\end{definition}

\begin{definition} \label{def:automaton}
A finite automaton is defined by
\begin{itemize}
\item
A set of states $S$ with a start state $s_0 \in S$.
\item
An alphabet which we fix to be $\{0,1\}$,
\item
A transition function $\delta : S \times \{0,1\} \to S$, $\delta(s,b)$ is the state
the automaton is, given that it was in state $s$ and the current input symbol is $b$.
For a string $w = w_1 \ldots w_k \in \{0,1\}^{\ast}$ we let $\delta(s,w)$ be defined
inductively by $\delta(s,w_1 \ldots w_k) = \delta(\delta(s, w_1), w_2 \ldots w_k)$.
\item
Two disjoint sets of final states $S_0$ and $S_1$.
The automaton stops whenever it is in $s \in S_0 \cup S_1$ (formally, for all $s \in S_0 \cup S_1$
and $b \in \{0,1\}$, it holds that $\delta(s,b) = s$).
If it stopped at $S_0$, then the output is $0$; if it stopped in $S_1$, then the output is $1$.
\end{itemize}
\end{definition}

This is a slight variation of the standard definition (see e.g. \cite{HU}) 
where the automaton stops when the it gets to the end of the input.
Here, the input is unbounded. 

Letting $L_i$ be the strings in $\{0,1\}^{\ast}$ 
for which the automaton stops at $S_i$, 
the automaton will define a simulation
of a function $f : (0,1) \to (0,1)$, if $\P_p[L_0 \cup L_1] = 1$.
Note that $L_0$ and $L_1$ are both regular languages, and $L_0 \cup L_1$ has
the prefix property. Moreover, for every pair of disjoint regular 
languages $L_0$ and $L_1$ with the prefix property, using
a standard product construction, it easy to write down an automaton as in Definition \ref{def:automaton} which
outputs $0$ for $L_0$ and $1$ for $L_1$.

Note that in Figure 1, for the leftmost automaton,
$L_0 = (00 + 11)^{\ast} 01$, $L_1 = (00 + 11)^{\ast} 10$; 
for the middle automaton,
$L_0 = 11 (0 + 1)^{\ast}$, 
$L_1 = (00 + 01 + 10) (0 + 1)^{\ast}$; 
and for the rightmost automaton,
$L_0 = (01 + 10)^{\ast} 00$, $L_1 = (01 + 10)^{\ast} 11$.

The definition below is a slight variation of the standard definition
of Pushdown automata (see e.g. \cite{HU} for the standard definition).
\begin{definition} \label{def:push}
A {\em Pushdown automaton} over the alphabet $\{0,1\}$ is defined by
\begin{itemize}
\item
A set of states $S$ with a start state $s_0 \in Q$.
\item
A stack alphabet $\Lambda$.
\item
A stack $T \in \Lambda^\ast$ which is initialized as the non-empty string $\tau$.
\item
A transition function $(\delta,\eta) : 
S \times \{0,1\} \times \Lambda  \to S \times \Lambda^{\ast}$.
The first coordinate, $\delta(s,b,\tau)$, is the new state
of the automaton given that it was in state $s$, the current input symbol is $b$ and that the
symbol at the top of the stack is $\tau$.
After this transition the symbol at the top of the stack is replaced
by the string in the second coordinate of $\eta(s,b,\tau)$.
\item
Two disjoint sets of final states $S_0$ and $S_1$. The automaton stops when the stack is empty.
If the automaton stops at a state of $S_0$, the output is $0$, and if
it stops at a state of $S_1$,
the output is $1$.
\end{itemize}
\end{definition}

Note that if $L_0$ ($L_1$) are the languages where the automaton output $0$ ($1$), then $L_0$ and $L_1$
are context free languages, and $L_0 \cup L_1$ has the prefix property.

\subsection{Motivation and related models}
Our results have a similar flavor to classical results of the Chomsky Sch\"{u}tzenberger theory \cite{CS} - thus a
relationship is established between the strength of a computation model and algebraic properties of the outputs
generated by this model. We discuss some of the connections in more
detail in subsections \ref{subsec:algauto} and
\ref{subsec:cs}.

The model we introduce in this paper has also some interesting relationships with the theory of computability
\cite{Gr,PR} as in both models the interest is in "real" inputs and outputs. In particular for both models,
the input is unbounded. See Subsection \ref{subsec:comp} for a more formal discussion.

Finally a major motivation for studying this model is the theory of
exact sampling, see e.g.
\cite{Bo,Al,PW1,PW2}. 
In the theory of exact sampling 
the aim is to sample (exactly) from a certain complicated
distribution given a simpler given distribution. 
Thus the problem is a basic problem in
the theory of exact sampling. In Subsection \ref{subsec:ex} we give a (toy) example which illustrating how our results
may be used for exact sampling.

\subsection{Paper plan}
In Section \ref{sec:auto} we prove Theorem \ref{thm:auto} and its generalizations to collections of coins, dice etc.
In Section \ref{sec:push} we prove Theorem \ref{thm:newpush} and 
show how using a push-down automaton, 
it is possible to simulate non-rational constants
and functions. In Section \ref{sec:related}
we discuss the relationships with the theory of languages and their formal power series, the theory of computability,
exact sampling and suggest some open problems.

\section{Rationality and finite automata} \label{sec:auto}

A sub-family of the coins which can be simulated via 
finite-automata are those which can be simulated via blocks.
Von Neumann's trick consists of reading $2$ bit blocks until 
$01$ or $10$ are reached, and then deciding
$0$ for $01$ and $1$ for $10$. 
Block simulation is a generalization of this procedure defined as follows
\begin{definition} \label{def:block_simuation}
A {\em block simulation} of $f$, is a simulation of $f$ of the following form.
Let $A_0$ and $A_1$ be disjoint subsets of $\{0,1\}^k$, and $A' = \{0,1\}^k \setminus (A_0 \cup A_1)$.
The simulation procedure has $L_0 = (\sum_{w \in A'} w)^{\ast} (\sum_{w \in A_0} w)$ and
$L_1 = (\sum_{w \in A'} w)^{\ast} (\sum_{w \in A_1} w)$.
In other words, the procedure reads a $k$ bit string $w$.
If $w \in A_0$, the procedure outputs $0$, if $w \in A_1$
the procedure outputs $1$; otherwise the procedure discards $w$ and reads a new $k$ bit string.
\end{definition}

Below we prove the following theorem which immediately implies 
Theorem \ref{thm:newauto}.

\begin{theorem} \label{thm:auto}
Let $\D \subset (0,1)$.
For $f : \D \to (0,1)$ the following are equivalent.
\begin{itemize}
\item[I]
$f$ can be block simulated.
\item[II]
$f$ can be simulated via a finite automaton.
\item[III]
$f$ is the restriction to $\D$ of a a
rational function $F$ over $\Q$ such that
$0 < F(p) < 1$ for all $p \in (0,1)$.
\end{itemize}
\end{theorem}

Note that $I \Rightarrow II$ is trivial. 
The implication $II \Rightarrow III$ is fairly easy.
We know of no way of proving $II \Rightarrow I$ directly. Instead, we prove $III \Rightarrow I$.
We know of no simple bound in terms of $f$ on the size of automaton (or block) needed in order to simulate $f$,
although there exists a simple algorithm for constructing an automaton that simulates $f$.

\subsection{Finite automaton $\Rightarrow$ Rationality}
\begin{proposition} \label{prop:auto}
Let $\D \subset (0,1)$. If a finite automaton $\Sigma$ simulates 
$f : \D \to (0,1)$, then $f$ is the restriction to $\D$ of a 
rational function $F$ over $\Q$ such that
$0 < F(p) < 1$ for all $p \in (0,1)$. 
Moreover if $\Sigma$ has $n$ states, then $F(p) = g(p)/h(p)$, where
$g(p),h(p) \in \Z[p]$ and $g$ and $h$ are of degree at most $n$.
\end{proposition}
The proof applies the maximum principle
for harmonic functions on directed graphs.
\begin{lemma} \label{lem:max}
Let $\Sigma$ be a finite automaton and $\FF$ a set of states such that for all $s \in S$, there exists
$w \in \{0,1\}^{\ast}$, such that $\delta(s,w) \in \FF$.
Let $0 < p < 1$ and $f : S \to \R$ be harmonic, so for all $s$ it holds that
$f(s) = p f (\delta(s,1)) + (1 - p) f(\delta(s,0))$.
Then $f$ achieves its maximum and minimum in $\FF$.
Moreover, given the values of $f$ in $\FF$, the values of $f$ in $S$ are uniquely determined.
\end{lemma}

\proofs
Note first that the last assertion follows from the first one, as if $f_1$ and $f_2$ are two
harmonic functions which have the same value on $\FF$, then $f_1 - f_2$ is harmonic and has the value
$0$ on $\FF$, which implies by the first assertion that $\max f_1 - f_2 = \min f_1 - f_2 = 0$, or $f_1 = f_2$.

In order to prove the first assertion, let $m = \max f$. 
Note that if $f(s) = m$, then $f(\delta(s,0)) = f(\delta(s,1)) = m$.
Letting $w \in \{0,1\}^{\ast}$ be such that $\delta(s,w) \in \FF$, we obtain
$f(\delta(s,w)) = m$ as needed.
$\QED$

\prooft{of Proposition \ref{prop:auto}}
Suppose that $f$ can is simulated by a finite automaton $\Sigma$. 
Let $S'$ be the set of states $s$ such that there exists
$w \in \{0,1\}^{\ast}$ with $\delta(s_0,w) = s$.
Clearly, we may remove from the automaton all the states not in $S'$
(redefining $\delta$ by restriction) and still obtain a finite automaton which simulates $f(p)$.
By the assumption that $\P_p[L_0 \cup L_1] = 1$, it follows that for all $s \in S'$, there
exists $w \in \{0,1\}^{\ast}$ such that $\delta(s,w) \in S_0 \cup
S_1$. From now on we assume that $S = S'$.

Note that since for all $s \in S$, there exists $w \in
\{0,1\}^{\ast}$ such that $\delta(s,w) \in S_0 \cup S_1$, it follows
that $\P_p[L_0 \cup L_1] = 1$ for all $0 < p < 1$.  Therefore the function $f$
is the restriction of the function $F(p) = \P_p[L_1]$ to $\D$. 

Suppose that $F(s;p)$ satisfies $F(s;p) = 0$ for $s \in S_0$, 
and $F(s;p) = 1$ for $s \in S_1$. For all other $s$, assume that
$$
F(s;p) = p F(\delta(s,1);p) + (1 - p) F(\delta(s,0);p) \,.
$$ 
By Lemma \ref{lem:max}, the function $F(s;p)$ 
is uniquely determined by these equations.
This implies that $F(s_0;p) = F(p)$.
Since $F(p)$ is uniquely determined by a collection of 
linear equation with coefficients in $\Z[p]$,
it follows by Cramer's rule that $F(p)$ may be written as the ratio of two determinants in $\Z[p]$,
and therefore $F(p) = g(p)/h(p) \in \Q(p)$, where the degrees of $g$ and $h$ are at most the number of states
of the automaton, as needed.
$\QED$

\subsection{Block Simulation}
In this subsection we study what can be simulated by blocks.

\begin{proposition} \label{prop:block}
$f$ can be simulated using a block procedure if and only if
$f$ can be written as $D(p)/E(p)$ where
\begin{eqnarray} \label{eq:DE}
D(p) &=& \sum_{i=0}^k d_i p^i (1-p)^{k-i}, \\ \nonumber
E(p) &=& \sum_{i=0}^k e_i p^i (1-p)^{k-i},
\end{eqnarray}
and for all $i$, the coefficients $d_i$ and $e_i$ are integers such that
$0 \leq d_i \leq e_i$.
\end{proposition}


\proofs
Suppose that $f$ is block simulated.
For a string $w \in \{0,1\}^{\ast}$, we write $n_1(w)$ for the number of $1$s in $w$.
Then
\[
f(p) = \P_p[L_1] = \frac{\sum_{w \in A_1} \P_p[w]}{\sum_{w \in A_0 \cup A_1} \P_p[w]} =
     \frac{\sum_{i=0}^k d_i p^i (1-p)^{k-i}}{\sum_{i=0}^k e_i p^i (1-p)^{k-i}},
\]
where
\[
d_i = \# \{w \in A_1 : n_1(w) = i\},
\]
and
\[
e_i = \# \{w \in A_1 \cup A_0 : n_1(w) = i\},
\]
satisfy that $0 \leq d_i \leq e_i$, as needed.

For the other direction, suppose that $f(p) = D(p)/E(p)$, where $D$ and $E$ satisfy (\ref{eq:DE}).
Let $r$ be a number such that $e_i \leq \binom{k}{i} \binom{2r}{r}$
for all $i$.
For each $i$, fix a bijection $B_i$ from
\[
 \{w \in \{0,1\}^k : n_1(w) = i\} \times \{v \in \{0,1\}^{2r} : n_1(v) = r\}
\]
to $\{1,\ldots,\binom{k}{i} \binom{2r}{r}\}$.

The sets $A_0$ and $A_1$ are subsets of $\{0,1\}^k \times \{0,1\}^{2r}$,
defined as follows. $A_1$ is defined as
\[ 
\cup_{i=0}^k 
\{(v,w) : n_1(v) = i, n_1(w) = r,  \mbox{ and } B_i(v,w) \leq d_i \}
\]
and $A_0$ as
\[
\cup_{i=0}^k
\{(v,w) : n_1(v) = i, n_1(w) = r, \mbox{ and } d_i < B_i(v,w) \leq e_i\}.
\]

So
\begin{eqnarray*}
\P_p[A_1] &=& p^r (1-p)^{r} \sum_{i=0}^k d_i p^i (1-p)^{k-i}, \\
\P_p[A_0 \cup A_1] &=& p^r (1-p)^r \sum_{i=0}^k e_i p^i (1-p)^{k-i},
\end{eqnarray*}
and therefore,
\[
\P[L_1] = \frac{\P_p[A_1]}{\P_p[A_0 \cup A_1]} = f(p),
\]
as needed.
$\QED$

\subsection{Rationality $\Rightarrow$ Finite automaton}
In this section we prove Theorem \ref{thm:auto}. The proof is based on a beautiful theorem by P\'{o}lya \cite{Po}
(see \cite{HLP}, 57--59).
We let $\Delta^s$ denote the open $s$-simplex of probability distributions,
\[
\Delta^s = \{p \in (0,1)^{s+1} : \sum_{i=1}^{s+1} p_i = 1\}.
\]

\begin{theorem}[P\'{o}lya \cite{Po}] \label{thm:polya}
Let $f : \Delta^{s-1} \to \R$ be homogeneous and positive polynomial in the variables $p_1,\ldots,p_s$.
Then for all sufficiently large $n$, all the coefficients of $(p_1 + \ldots + p_s)^n f(p_1,\ldots,p_s)$ are
positive.
\end{theorem}

\begin{lemma} \label{lem:posrat}
Let $f : (0,1) \to (0,1)$ be a rational function. Then there exist polynomials $d$ and $e$
\begin{eqnarray} \label{eq:DR2}
d(p) &=& \sum_{i=0}^k d_i p^i (1-p)^{k-i}, \\ \nonumber
e(p) &=& \sum_{i=0}^k e_i p^i (1-p)^{k-i},
\end{eqnarray}
where for all $i$, the coefficients $d_i$ and $e_i$ are integers such that
$0 \leq d_i \leq e_i$, and $f(p) = d(p)/e(p)$.
\end{lemma}

\proofs
As $f(p)$ is a rational function it may be written in the form 
$\overline{D}(p)/\overline{E}(p)$, where $\overline{D}(p) \in \Z[p]$ and
$\overline{E}(p) \in \Z[p]$ are relatively prime polynomials. 
Since $0 < f(p)$ for all $0 < p < 1$, 
it follows that $\overline{D}(p)$ and $\overline{E}(p)$ do not change sign
in the interval $(0,1)$. 
Without loss of generality we assume that $\overline{D}(p) > 0$ 
and $\overline{E}(p) > 0$ for all $p \in (0,1)$. 
Note furthermore that if $\overline{D}(p) = \sum_{i=0}^k a_i p^i$
and $\overline{E}(p) = \sum_{i=0}^k b_i p^i$,
then we may define homogeneous polynomials $D(p,q)$ and $E(p,q)$ 
of degree $k$, by letting 
$D(p,q) = \sum_{i=0}^k a_i p^i (p + q)^{k-i}$, and 
$E(p,q) = \sum_{i=0}^k b_i p^i (p + q)^{k-i}$.
Note that $\overline{D}(p) = D(p,1-p)$ and $\overline{E}(p) = E(p,1-p)$.
Let us rewrite,
\begin{eqnarray}
D(p,q) &=& \sum_{i=0}^k d_i p^i q^{k-i}, \\ \nonumber
E(p,q) &=& \sum_{i=0}^k e_i p^i q^{k-i}.
\end{eqnarray}
The polynomials $D(p,q), E(p,q)$ and $E(p,q) - D(p,q)$ 
are all positive homogeneous polynomials.
Therefore by Theorem \ref{thm:polya}, if follows that there exists an $n$ such that
letting $d(p,q) = (p + q)^n D(p,q)$ and $e(p,q) = (p + q)^n E(p,q)$, 
the polynomials $d$,$e$ and $e-d$
all have positive coefficients as polynomials in $p$ and $q$.
Writing $f(p) = d(p,1-p)/e(p,1-p)$ we obtain the required result.
$\QED$

\prooft{of Theorem \ref{thm:auto}}
The implication $I \Rightarrow II$ is trivial, the implication $II \Rightarrow III$ follows
from Proposition \ref{prop:auto}, while $III \Rightarrow I$ follows from Lemma \ref{lem:posrat} together
with Proposition \ref{prop:block}.

\subsection{Extensions to dice and other $k$-sided coins}
In this subsection we discuss generalizations of the problem to $k$-sided coins, such as dice.
A {\em simulation} of
$f = (f_1,\ldots,f_t) : \Delta^s \to \Delta^t$, is a collection of $t$ disjoint languages
property $(L_1,\ldots,L_t)$ over the alphabet $\Sigma =
\{1,\ldots,s\}$,  such that
$\cup_{i=1}^t L_i$ has the prefix property and
$\P_p[L_i] = f_i(p)$ for all $p \in \Delta^s$. 
The definition of simulation via finite/pushdown automata and
Turing machines naturally extend to this setting. The continuity results of \cite{KO} extend to the more general
setting as well.

\begin{proposition} \label{prop:auto2}
If a finite automaton simulates $f : \Delta^s \to \Delta^t$,
then $f$ is a rational function over $\Z$
(i.e. $f_i(p)$ is a rational function
over $\Z$ for all $1 \leq i \leq t+1$).
\end{proposition}
\proofs
Repeat the proof of Proposition \ref{prop:auto} for each of the $f_i$'s.
$\QED$.

Repeating the proof of Theorem \ref{thm:auto}, we see that
\begin{theorem} \label{thm:auto2}
Any rational function $f : \Delta^s \to \Delta^1$ can be simulated via blocks.
\end{theorem}

From which we conclude that
\begin{corollary} \label{cor:auto2}
Any rational function $f : \Delta^s \to \Delta^t$ can be simulated via blocks.
\end{corollary}

\begin{proof}
The proof is by induction on $t$.
The case $t=1$ is covered by Theorem \ref{thm:auto2}.
Suppose $t > 1$, and let $(f_1,\ldots,f_{t+1}) : \Delta^s \to \Delta^t$ be a rational function.
By Theorem \ref{thm:auto2}, there exists a block simulation $(A_1,A_2) \subset (\Sigma^k)^2$ for $(f_1,1 - f_1)$,
and by the induction hypothesis there exists a block simulation $(B_1,\ldots,B_t) \subset (\Sigma^r)^t$
for $(f_2/(1 - f_1),\ldots,f_{t+1}/(1 - f_1))$. Taking
\[
(A_1 \times \Sigma^r, A_2 \times B_1, A_2 \times B_2, \ldots, A_2 \times B_t) \subset (\Sigma^{r+k})^{t+1},
\]
we obtain a block simulation for $(f_1,\ldots,f_{t+1})$ as needed.
\end{proof}

\section{Pushdown automata} \label{sec:push}

We now prove Theorem \ref{thm:newpush}.   
We begin by showing that if $f$ is simulated by a
pushdown automaton, then $f$ is the
unique solution of a set of polynomial equations. We then invoke the
results of the appendix to deduce that $f$ is an algebraic function.
In Section \ref{sec:push} we also construct pushdown automata which
simulate non-rational functions such as $f(p) = \sqrt{p}$.
We don't know if every algebraic $f : (0,1) \to (0,1)$ can be 
simulated by a pushdown automaton.



\subsection{Pushdown automata and algebraic functions} \label{subsec:algauto}

The Chomsky-Sch\"{u}tzenberger theory implies that if $L_0$ and $L_1$ are languages which are generated by unambiguous
grammars and $f(p) = \P_p[L_1]$, then $f(p)$ is an algebraic function. In this subsection we aim to prove algebraic properties
of $f$ even when $L_0$ and $L_1$ are inherently ambiguous.

Suppose $\Sigma$ is a pushdown automaton which simulates a function $f$.
Call $(b,s) \in \Lambda \times S$ {\em good}, if when the automaton is at state $s$ and the
stack is $bw$ (where $b$ is at the top), then with probability $1$ at some point the stack will be $w$.
Call $(b,s)$ {\em bad} otherwise.
By the assumption that $\P_p[L_0 \cup L_1] = 1$, it follows that starting at $(s_0,\tau)$ it is impossible
for the automaton to reach a state $s$ with $bw$ at the top of the stack, where $(b,s)$ is bad.
Thus we can redefine all transitions $(b,s) \to (b'w',s')$, where $(b',s')$ is bad, in an arbitrary manner,
and still obtain $\P_p[L_0 \cup L_1] = 1$. Therefore, without loss of generality we may assume that all
$(b,s) \in \Lambda \times S$ are good.

Let 
$\alpha(p ; b,s,s'),\; \alpha : [0,1] \times \Lambda \times S \times S
\to [0,1]$ 
be defined as follows.
For $w \in \Lambda^{\ast}$, let $\alpha(p ; b,s,s')$ 
be the probability that given that
currently the automaton is at state $s$ and has in its stack $bw$ 
(where $b$ is at the top),
at the first time that the content of the stack will be $w$, 
it will be at state $s'$.
It is easily seen that $\alpha(p ; b,s,s')$ is well defined (does not depend on $w$). Moreover, by the assumption
that all $(b,s)$ are good, it follows that 
$\sum_{s' \in S_0 \cup S_1} \alpha(p ; b,s,s') = 1$.

We extend the definition of $\alpha$ to 
$\tilde{\alpha}(p ; u,s,s'),\; \tilde{\alpha} : \Lambda^{\ast} \times S \times S \to [0,1]$,
where $\tilde{\alpha}(p ; u,s,s')$ is the probability that given that
currently the automaton is at state $s$ and has in its stack $u w$ 
(where $u$ is above $w$),
at the first time that the content of the stack will be $w$, 
it will be at state $s'$.
Note that if $w = w_1 \ldots w_r$, then
\begin{equation} \label{eq:defalph}
\tilde{\alpha}(p ; w,s,s') = \sum
                     \prod_{i=1}^r \alpha(p ; w_i,s_{i},s_{i+1}),
\end{equation}
where the sum is over all
\[
(s_1,s_2,\ldots,s_r,s_{r+1}) \in \{s\} \times S^{r-1} \times \{s'\},
\]
and if $\eps$ denotes the empty word, then
\begin{equation} \label{eq:defalph2}
\forall s,s' \in Q:\, 
\tilde{\alpha}(p; \epsilon,s,s') = 1_{\{s = s'\}}.
\end{equation} 

Note that if $\tau$ is the initial word at the stack, then
\begin{equation} \label{eq:l1stack}
\P_p[L_1] = \sum_{s' \in S_1} \tilde{\alpha}(p ; \tau, s_0, s').
\end{equation}

Therefore if we could prove algebraic properties of the functions 
$\alpha(p ; b,s,s')$,
we will deduce algebraic properties of $f$.

\begin{claim} \label{claim:inf_harmonic}
For all $0 < p < 1$,
$\tilde{\alpha}(p; \cdot, \cdot, \cdot) : 
\Lambda^{\ast} \times S \times S \to \R$ is the unique
bounded solution of the equations
\begin{eqnarray} \label{eq:alpha_harmonic}
\forall w \in \Lambda^{\ast},\, \forall b \in \Lambda ,\, \forall s,s' 
\in Q:\, 
\tilde{\alpha}(p; bw, s, s') &=&  
p      \tilde{\alpha} \left(p; \delta(s,1,b) w,\eta(s,1,b),s'\right) 
\\ \nonumber &+&
(1-p)  \tilde{\alpha} \left(p; \delta(s,0,b) w,\eta(s,0,b),s'\right),
\end{eqnarray}
\begin{equation} \label{eq:alpha_boundary}
\forall s,s' \in Q:\, 
\tilde{\alpha}(p; \epsilon,s,s') = 1_{\{s = s'\}} \,\,\, 
(\epsilon \mbox{ is the empty word}).
\end{equation}
\end{claim}

\proofs
By linearity, it suffices to prove that the zero function is the only bounded 
solution to (\ref{eq:alpha_harmonic}) with the boundary conditions 
\begin{equation} \label{eq:alpha_boundary_0}
\forall s,s' \in Q:\, 
\tilde{\alpha}(p; \epsilon,s,s') = 0.
\end{equation}
Fix $(w,s) \in \Lambda^{\ast} \times Q$. We will show that 
$\tilde{\alpha}(p; w,s,s') = 0$.
Consider the random walk $(W_t,S_t)_{w,s}$ defined on the graph 
$\Lambda^{\ast} \times Q$, where $(W_0,S_0)_{w,s} = (w,s)$.  
Given $W_t = B U$, where $B \in \Lambda$, 
the conditional probabilities for $(W_{t+1},S_{t+1})$ are given by
\begin{equation} \label{eq:steprw}
(W_{t+1},S_{t+1}) = \left\{ \begin{array}{ll}
                    \left(\delta(S_t,1,B) U, \eta(S_t,1,B)
                    \right) & \mbox{ with probability } p, \\
                    \left(\delta(S_t,0,B) U, \eta(S_t,0,B)
                    \right) & \mbox{ with probability } 1-p.
                    \end{array} \right. 
\end{equation}
If $W_t$ is the empty word, then we let $(W_{t+1},S_{t+1}) =
(W_t,S_t)$.
By definition, for all $s' \in Q$, the process $\tilde{\alpha}(p; W_t,S_t,s')$
is a bounded martingale. The assumption that the pushdown automaton
stops a.s. implies by (\ref{eq:alpha_boundary_0}) that the martingale
converges to $0$ a.s. (and therefore in $L_1$). We therefore conclude
that $\tilde{\alpha}(p; W_t, S_t, s')$ is identically $0$ as needed.
$\QED$

Recall that given that the current state is $s$, 
and the top of the stack is $b$, with probability
$p$ the automaton will move to state $s_1 = \delta(s,1,b)$, and
instead of $b$, the top of the stack will contain
$\eta(s,1,b) = c_1 \ldots c_r$;
with probability $1-p$ the automaton will move to state
$\bar{s}_1 = \delta(s,0,b)$, 
and instead of $b$, the top of the stack will contain
$\eta(s,0,b) = {\bar{c}}_1 \ldots {\bar{c}}_{\bar{r}}$.

We can therefore write
\begin{eqnarray} \label{eq:alph_recur1}
\alpha(p ; b,s,s') &=&
p  \left( 1_{\{r = 0\}} 1_{\{s_1 = s'\}} + 1_{\{r > 0\}} 
   \tilde{\alpha}(p ; c_1 \ldots c_r,s_1,s') \right) \\ \nonumber &+&
(1-p) \left( 
1_{\{\bar{r} = 0\}} 1_{\{\bar{s}_1 = s'\}} + 1_{\{\bar{r} > 0\}} 
\tilde{\alpha}(p ; \bar{c}_1 \ldots \bar{c}_{\bar{r}},s_1,s') \right),
\end{eqnarray}
Or
\begin{eqnarray} \label{eq:alph_recur2}
\alpha(p ; b,s,s') &=&
p \left( 1_{\{r = 0\}} 1_{\{s_1 = s'\}} + 1_{\{r > 0\}} \sum
\prod_{i=1}^r \alpha(p ; c_i,s_{i},s_{i+1}) \right) \\ \nonumber &+&
(1-p) \left( 
1_{\{\bar{r} = 0\}} 1_{\{\bar{s}_1 = s'\}} + 1_{\{\bar{r} > 0\}} 
\sum \prod_{i=1}^{\bar{r}} \alpha(p ;
{\bar{c}}_i,\bar{s}_{i},\bar{s}_{i+1})
\right).
\end{eqnarray}
where the first (second) sum is taken over all
\begin{eqnarray*}
(s_1,\ldots,s_{r+1}) &\in& \{s_1\} \times S^{r-1} \times \{s'\}, \\
(\bar{s}_1,\ldots,\bar{s}_{r'+1}) &\in& 
\{\bar{s}_1\} \times S^{\bar{r}-1} \times \{s'\}.
\end{eqnarray*}
Note that (\ref{eq:alph_recur2}) defines a set of algebraic equations
in $p$ and $\{\alpha(p ; b,s,s')\}_{b,s,s'}$.

\begin{claim} \label{claim:eq_alpha}
For all $0 < p < 1$, there is a unique positive solution
$\alpha(p; \cdot, \cdot, \cdot) : [0,1] \times \Lambda \times S \times
S \to \R$, to equations (\ref{eq:alph_recur2}) and 
\begin{equation} \label{eq:alpha_prob}
\forall b \in \Lambda,\,\forall s \in Q:\,
\sum_{s' \in Q} \alpha(p;b,s,s') = 1.
\end{equation}
\end{claim}

\proofs
In order to prove the claim it suffices to show that each positive
solution to 
(\ref{eq:alph_recur2}) and (\ref{eq:alpha_prob}) defines a positive
bounded solution to (\ref{eq:alpha_harmonic}) and
(\ref{eq:alpha_boundary}) via (\ref{eq:defalph}) and
(\ref{eq:defalph2}). 

It is immediate to see that (\ref{eq:alph_recur2}) implies 
(\ref{eq:alpha_harmonic}). Moreover, 
(\ref{eq:alpha_prob}) implies by (\ref{eq:defalph}) that for all 
$w \in \Lambda^{\ast}$ and $s \in Q$, it holds that 
\[
\sum_{s' \in Q} \tilde{\alpha}(p;w,s,s') = 1.
\]
and so $\tilde{\alpha}$ is a bounded function as needed.
$\QED$

\prooft{of Theorem \ref{thm:newpush}}
The result follows immediately from Claim \ref{claim:eq_alpha} by 
Theorem \ref{thm1}.
$\QED$

\subsection{Pushdown automata which simulate non-rational functions} \label{subsec:pushex}
In this subsection we construct a pushdown automaton which simulates a non-rational function.
Let $g : (0,1) \to (0,1/2)$ be a rational function. We'll construct a pushdown automata which simulates
the function $\gamma(p) = \frac{1 - \sqrt{ 1 - 2g(p)}}{4 g(p)}$.

Taking $g(p) = (1 - p)/2$,
we obtain $\gamma(p) = (1 - \sqrt{p})/(1 - p)$ - thus using a product
construction, it is easy to construct a pushdown automaton simulating
the function $f(p) = \sqrt{p}$. 

\begin{center}
\begin{picture}(5,5)
\thicklines
\put(1,0.5){$0$}
\put(2,0.5){$1$}
\put(1,1){\line(1,0){1}}
\put(1,1.5){\line(1,0){1}}
\put(1,2){\line(1,0){1}}
\put(1,2.5){\line(1,0){1}}
\put(1,3){\line(1,0){1}}
\put(1,3.5){\line(1,0){1}}
\put(1,1){\line(0,1){3}}
\put(2,1){\line(0,1){3}}
\put(3,2){\vector(1,0){0.5}}
\put(3,2){\vector(0,1){0.5}}
\put(3,3){$g(p)$}
\put(3,1){$g(p)$}
\put(5,1){$g(p)$}
\put(5,3){$g(p)$}
\put(3.25,2.2){$1 - 2 g(p)$}
\put(3,2){\vector(0,-1){0.5}}
\put(5,2){\vector(-1,0){0.5}}
\put(5,2){\vector(0,1){0.5}}
\put(5,2){\vector(0,-1){0.5}}
\end{picture}
\end{center}

Consider a random walk on the ladder graph $\N \times \{0,1\}$ where an edge $((x,y),(x',y'))$ is present if
$x = x'$ and $|y - y'| = 1$, or $|x - x'| = 1$ and $y = y'$. The random walk moves to the left (right) with
probability $1 - 2 g(p)$, and up (down) with probability $g(p)$.
Let $\gamma$ be the probability that starting at $(0,1)$
the first hitting point of the random walk at level $0$ ($\{(0,0),(1,0)\}$) is $(0,0)$.
It is easy to see that
\begin{eqnarray*}
\gamma &=& g(p) + (1 - 2 g(p)) (1 - \gamma) + g(p) 
(\gamma^2 + (1 - \gamma)^2) \\ &=& 1 - \gamma + 2 g(p) \gamma^2,
\end{eqnarray*}
and therefore
\begin{equation} \label{eq:biasq}
\gamma(p) = \frac{1 - \sqrt{ 1 - 2 g(p)}}{2 g(p)}.
\end{equation}

It is easy to simulate the random walk with a pushdown automaton.
The stack alphabet is $\{x\}$ where $x^n$ correspond to level $n$ of the ladder -- the initial word
at the stack is $x$.

Assume first that the input alphabet is $\{0,1,2\}$ where the probability that letter $1$ appears is $1 - 2 g(p)$,
and the probability that letter $0$ or $2$ appears is $g(p)$.

Let $s_0$ and $s_1$ be two states of the automaton corresponding the left
$(\{0\} \times \N)$ and right $(\{1\} \times \N)$ of the ladder.
Reading the symbol $1$ will correspond
to a transition from $s_0$ to $s_1$ or vice-versa without changing the content of the stack.
Reading the symbol $2$ at state $s_i$ will result at staying at state $s_i$ and pushing an $x$ to the stack.
Reading the symbol $0$ at state $s_i$ with $x$ at the top of the stack, will result at staying
at state $s_i$ while popping $x$ from the stack.
In this way it is possible to simulate the random walk given an infinite sequence of $\{0,1,2\}$ symbols with
$(g(p), 1 - 2g(p), g(p))$ bias - and therefore toss a coin with bias (\ref{eq:biasq}).

In the general case where we are given an infinite sequence of $(p,1-p)$ bits, we use block constructions
of Section \ref{sec:auto} in order to
generate a sequence of $(g(p), 1 - 2g(p), g(p))$ $\{0,1,2\}$ variables together with the above construction in order
to obtain the required result.

\begin{remark}
Similarly one may construct a pushdown automaton associated with
a random walk on the ladder where the probabilities of going
$(\mbox{up,left-right,down})$ are given by $(g(p), 1 - g(p) -
h(p),h(p))$. Note however that in this case $\P_p[L_0 \cup L_1] = 1$ 
{\bf iff} the random walk is recurrent {\bf iff} $h(p) \geq g(p)$.
Thus, unlike finite automata, there exist 
pushdown automata which define a valid simulation
only for a proper subset of the interval $(0,1)$.
\end{remark}   

\section{Related models and open problems} \label{sec:related}

\subsection{Chomsky-Sch\"{u}tzenberger theory} \label{subsec:cs}

In the seminal paper \cite{CS} titled 
``The algebraic theory of context-free languages''
by Chomsky and Sch\"{u}tzenberger,
the authors discuss many beautiful relationships between different types of languages and their generating functions.
The generating function of a grammar generating a language $L$ is defined as $f = \sum_{w \in L} n(w) w$, where $n(w)$
is the number of derivations of the word $w$ in the grammar, 
and $f$ is viewed as a formal power series in
the non-commutative variables $0$ and $1$. 
The results of \cite{CS} imply in particular that the
generating function of a regular language is rational, 
which implies Proposition \ref{prop:auto} --
simulation via finite automata always yield rational functions.

In \cite{CS} it is also proven that the generating function of a 
context-free language is algebraic.
However, this does not imply that coins tossed via pushdown automata 
are algebraic, as many of the context free
languages are inherently ambiguous and for such languages, the 
non-commutative power series 
$P(L) = \sum_{w \in \{0,1\}^{\ast}} 1_{\{w \in L\}} w$ is not algebraic.
(see e.g. \cite{HU}; as Larry Ruzzo kindly noted, there are also context free
languages with the prefix property that are inherently ambiguous).

Thus,
 while Proposition \ref{prop:auto} could be obtained by 
projecting the Chomsky-Sch\"{u}tzenberger results
from the non-commutative setting to the commutative setting, an analogous result for
pushdown automata cannot be obtained in a similar way.
Theorem \ref{thm:auto} may be interpreted as a commutative inverse to the Chomsky-Sch\"{u}tzenberger result for regular
languages, where additional positivity restrictions are imposed.

\subsection{Theory of computability and Turing machines} \label{subsec:comp}

A distinctive feature of the computation models discussed in this paper is that they have an unbounded input,
unlike the classical Turing machine for which the input is of bounded length.

In \cite{Gr} (see \cite{PR} for background), a
computational model in analysis is introduced. This model has some common features with the model introduced here --
in particular, the input is given as an unbounded sequence.
However, the models are different as our model
has an input drawn according to an i.i.d. distribution. Moreover, in our model we are looking for a coin which has
{\em exactly} the correct distribution $f(p)$, while
in \cite{Gr} the aim is find a Turing machine which for an (unbounded) input $x$, and an error bound $\eps$
computes an approximation of the function $f(x)$ within error margin $\eps$.

There are, however, some striking similarities between the two models. For instance the continuity results of
\cite{KO} should be compared with \cite{Gr} where it is shown that computable functions are continuous.
Moreover it is easy to adapt the proofs of \cite{KO} in order to show that if $f : (0,1) \to (0,1)$ is
computable and polynomially bounded from $0$ and $1$ at $0$ and $1$, then it possible to simulate the function $f(p)$
via a Turing machine.


\subsection{Exact sampling} \label{subsec:ex}
The theory of exact sampling (see e.g. \cite{Bo,Al,PW1,PW2}) 
deals with simulating a complicated probability measure
using a simple one. Usually, the simple measure consists of a 
sequence of unbiased bits or a a sequence of uniform variables.

In this paper both probability measures of interest 
(the $p$ coin and the $f(p)$ coin) are simple. The difficulty here
is that $p$ is unknown. It is easy to estimate $p$, and therefore
$f(p)$.  However, to get a coin with an exact bias $f(p)$ is harder.

Yet, the problem studied here may be interpreted as a problem in exact
sampling.
Consider the following examples.
\begin{itemize}
%

\item
Suppose that some physical process produces percolation configuration
on a grid $\Gamma_1$, where the probability of an open edge is $p$.
We are interested in performing percolation on a grid $\Gamma_2$,
where we want that the probability that an edge $e$ is open in
$\Gamma_2$ to be equal to the probability that two
vertices at distance $2$ in $\Gamma_1$ are connected by a path of
length $2$. Our results allow to use samples of the configuration 
$\Gamma_1$ in order to produce samples for the process on $\Gamma_2$.

\item
Let $\{x_i\}_{i \ge 1}$ be i.i.d. bits with unknown mean $\alpha$.
 Suppose that we are given the products $\{y_i\}_{i \ge 1}$, where 
$y_i=x_{2 i}  x_{2 i - 1}$,  
and we want to simulate (one or more) i.i.d.\ bits $\{z_j\}$ with mean $\alpha$. 
This can be done using the pushdown automaton in 
Subsection \ref{subsec:pushex} that simulates the function $\sqrt{p}$.
\end{itemize} 

\subsection{Some open problems} \label{subsec:open}

\begin{problem}
Let $f : (0,1) \to (0,1)$ be a rational function. What is the smallest size of an automaton that simulates $f(p)$?
Is there an efficient algorithm for finding this automaton?
\end{problem}
This problem is potentially hard, as the size of the automaton probably depends on analytic and number theoretic
properties of $f$ (see \cite{PR2}).

\begin{problem}
Let $f : (0,1) \to (0,1)$ be an algebraic function. Can $f$ be simulated by a pushdown automaton?
\end{problem}

{\bf Acknowledgment:} We thank Jim Propp for suggesting the
problem. We are grateful to Omer Angel, Paul Beame and Larry Ruzzo
for helpful remarks,
and to Russ Lyons and Sergey Fomin for references.
Most of this work was done while the first author was a postdoctoral
researcher at Microsoft Research, and the second author was visiting there.

\newpage
\appendix
\section{Appendix on Algebraic functions}
\begin{center}
Christopher J. Hillar
\footnote{This work is supported under a National Science Foundation
  Graduate Research Fellowship.} \\ 
chillar@math.berkeley.edu \\
\end{center} 
The purpose of this note is to establish the following fact using
techniques from real algebraic geometry.  It will be a direct
corollary of the more general Theorem \ref{moregenthm} below.
\begin{theorem}\label{thm1}
Let $\{F_i(p,x_1,\ldots,x_n)\}_{i=1}^m \subset \mathbb
Q[p,x_1,\ldots,x_n]$, and let $S$ be the set in $\mathbb R^{n+1}$
defined by
\begin{equation*}
S =  \left\{ (p,x_1,\ldots,x_n) : \ 0 < p < 1 \ \wedge \ \bigwedge
\limits_{j=1}^n {0 \leq x_j \leq 1} \ \wedge \ \bigwedge
\limits_{i=1}^m F_i(p,x_1,\ldots,x_n) = 0 \right\}.
\end{equation*}
Suppose that for each $p \in (0,1)$, there exist a unique
$(a_1,\ldots,a_n) \in \mathbb R^n$ such that $(p,a_1,\ldots,a_n)
\in S$; equivalently, $S$ is given as the image of some function
$\phi: (0,1) \rightarrow [0,1]^{n+1}$ with $\phi(p) =
(p,\phi_1(p),\ldots,\phi_n(p))$. Then, there exist nonzero
polynomials $g_j \in \mathbb Q[X,Y]$ such that
\[g_j(p,\phi_j(p)) = 0\] for all $p \in (0,1)$.
\end{theorem}

We begin with an abstract setting.  Let $K$ be a field, and let
$P$ be a subset of $K$ satisfying the following two properties:

\begin{enumerate}
  \item If $x,y \in P$, then $x+y$ and $xy \in P$.
  \item $K$ is the disjoint union of $P$, $\{0\}$, and $-P$.
\end{enumerate}
A subset $P$ as above is called the set of \emph{positive}
elements, and we say that $P$ is an \emph{ordering} of $K$.  A
field $K$ is then an \emph{ordered} field if there exists an
ordering of $K$.  A \emph{real field} is a field in which $-1$ is
not a sum of squares, and a \emph{real closed field} $R$ is a real
field such that such that any algebraic extension of $R$ that is
real must be equal to $R$.  For example, both $\mathbb R$ and
$\mathbb R \cap \mathbb Q^a$ are real closed fields (where
$\mathbb Q^a$ is the algebraic closure of $\mathbb Q$).  For ease
of notation below, we set $\mathbb Q^* = \mathbb R \cap \mathbb
Q^a$. From the definition, it is clear that any real field must
have characteristic 0, and thus $\mathbb Q$ naturally embeds in
any real field.  It also follows that $\mathbb Q^*$ is a subset of
every real closed field.

Any real closed field $R$ has a unique ordering, and the positive
elements are the squares of $R$.  Moreover, every polynomial of
odd degree in $R[x]$ has a root in $R$ \cite[p. 452]{Lang}.  In
light of these observations, the axioms for the theory of real
closed fields (RCF) consist of \cite[p. 24]{Haskell}:

\begin{enumerate}
  \item the axioms for ordered fields;
  \item $\forall x>0 \ \exists y \ y^2 = x$;
  \item the axiom $\forall x_0 \ldots \forall x_{n-1} \
        \exists y \ y^n + x_{n-1}y^{n-1} + \cdots +x_0 = 0$ for each odd $n > 0$.
\end{enumerate}

The set $S$ as in Theorem \ref{thm1} is an example of a set
definable by a Boolean combination of a finite number of
polynomial inequalities and equalities (over $\mathbb Q$).  Such a
set is called \emph{semialgebraic} (over $\mathbb Q$), and a
function is called \emph{semialgebraic} if its graph is a
semialgebraic set. We may now state the main theorem.

\begin{theorem}\label{moregenthm}
Let $R$ be a real closed field, and let $A$ be a semialgebraic
subset of $R$ defined using polynomials
$\left\{H_i(p)\right\}_{i=1}^k \subset \mathbb Q[p]$.  Also, let
$S$ be a semialgebraic subset of $R^{n+1}$ defined by
$\{F_i(p,x_1,\ldots,x_n)\}_{i=1}^m \subset \mathbb
Q[p,x_1,\ldots,x_n]$.  Suppose that for each $p \in A$, there
exist a unique point $(a_1,\ldots,a_n) \in R^n$ such that
$(p,a_1,\ldots,a_n) \in S$; equivalently, $S$ is given as the
image of some function $\phi: A \rightarrow R^{n+1}$ with $\phi(p)
= (p,\phi_1(p),\ldots,\phi_n(p))$.  Then, there exist nonzero
polynomials $g_j \in \mathbb Q[X,Y]$ such that
\[g_j(p,\phi_j(p)) = 0\] for all $p \in A$.
\end{theorem}

A fundamental fact about RCF is that it is a \emph{complete
theory} \cite[p. 19]{Haskell} in the sense that each first order
sentence expressible in the theory of RCF is either true in every
structure satisfying the RCF axioms or false in every such
structure. For example, the sentence, $\forall x \ \exists y \ y^2
= -x$, evaluates to false for any structure (such as $\mathbb R$)
satisfying the axioms of RCF.  As a standard application of
completeness, we present the following

\begin{lemma}\label{lem1}
Assume the hypothesis of Theorem \ref{moregenthm}.  Then, for each
$p \in A \cap \mathbb Q^*$, we have that $\phi_j(p) \in \mathbb
Q^*$ for all $j \in \{1,\ldots,n\}$.
\end{lemma}

\proofs To simplify notation, we write ``$p \in A$'', for example,
in place of the Boolean combination of polynomial inequalities and
equalities that defines $A$. As $\mathbb Q$ embeds in any real
closed field, the sentence,
\[\forall p \ \exists x_1 \ldots \exists  x_{n} \ \neg(p \in A) \ \vee \
(p,x_1,\ldots,x_n) \in S,\] is a valid sentence in any structure
satisfying the axioms of RCF. By completeness, it must have the
same truth value in every real closed field. Since it is a true
statement in $R$ by assumption, it follows that it is also true in
the real closed field $\mathbb Q^*$.  Let $p \in A \cap \mathbb
Q^*$. Then, there exists a tuple, $\mathbf{a} = (a_1,\ldots,a_n)
\in (\mathbb Q^*)^n$, such that $(p,a_1,\ldots,a_n) \in S$.  By
the hypothesis, this $\mathbf{a}$ must be the unique tuple in
$R^n$ for this $p$, and hence for all $j$, $\phi_j(p) = a_j \in
\mathbb Q^*$, completing the proof. $\QED$

A key result in the theory of RCF is the following theorem
(which is essentially a restatement of the fact that RCF has
quantifier elimination) \cite[p. 92]{Marker}.

\begin{theorem}[Tarski-Seidenberg Theorem]\label{thm2}
The projection of a semialgebraic set is semialgebraic.
\end{theorem}

\begin{corollary}\label{cor1}
The functions $\phi_j$ are semialgebraic functions.
\end{corollary}

\proofs Let $T = \left\{(p,x_1,\ldots,x_n) \ : \ p \in A \ \wedge
\ (p,x_1,\ldots,x_m) \in S \right\}$, which is semialgebraic.
Then, the image of the projection of $T$ into $R^2$ given by
$(p,a_1,\ldots,a_n) \mapsto (p,a_j)$ is also semialgebraic by the
theorem. $\QED$

Semialgebraic functions are well-behaved in the following sense
\cite[p. 17]{Marker2}.

\begin{theorem}\label{algthm}
Let $R$ be a real closed field.  If $T$ is a semialgebraic subset
of $R^n$ and $h : T \rightarrow R$ is semialgebraic, then there is
a nonzero polynomial $g(X_1,\ldots,X_n,Y) \in R[X_1,\ldots,X_n,Y]$
such that $g(\mathbf{x}, h(x)) = 0$ for all $\mathbf{x} \in T$.
\end{theorem}

We remark that applying Theorem \ref{algthm} with Corollary
\ref{cor1} already gives a result similar to Theorem
\ref{moregenthm}. The only subtlety is that we would like the
polynomial $g$ to be in $\mathbb Q[X,Y]$ instead of $R[X,Y]$. We
are now ready prove Theorem \ref{moregenthm}.

\proofs[Proof of Theorem \ref{moregenthm}] Fix $j \in
\{1,\ldots,n\}$.  We will apply Theorem \ref{algthm} with $R =
\mathbb Q^*$ and $T = A \cap \mathbb Q^*$. From Lemma \ref{lem1},
it follows that $\phi_j(T) \subseteq \mathbb Q^*$, and from
Corollary \ref{cor1}, we have that $\phi_j$ is semialgebraic.
Therefore, from Theorem \ref{algthm}, there is a nonzero
polynomial $g(X,Y) \in \mathbb Q^*[X,Y]$ such that $g(p,
\phi_j(p)) = 0$ for all $p \in A \cap \mathbb Q^*$.  We will now
produce a nonzero polynomial $\tilde{g} \in \mathbb Q[X,Y]$ with
the same property.

Consider the field, $\mathbb Q^*(X)$, of rational functions in the
variable $X$. View $g(X,Y)$ as a polynomial in the variable $Y$
over $\mathbb Q^*(X)$, and, upon clearing denominators, let
$q_i(X,Y) \in \mathbb Q^*[X,Y]$ ($i = 1,\ldots,r$) be the
irreducible factors of $g(X,Y)$ (over $\mathbb Q^*(X)$).  It is
clear that for all $p \in A \cap \mathbb Q^*$, we have
\[\prod\limits_{i = 1}^r {q_i(p,\phi_j(p))} = 0.\]  Fix $i$ and let
$\alpha_{1},\ldots,\alpha_t \in \mathbb Q^*$ be all of the
coefficients in $q_i(X,Y)$.  Extend $\mathbb Q(X)$ by these
coefficients, so that $\mathbb Q(X)(\alpha_1,\ldots,\alpha_t)$ is
a finite extension of $\mathbb Q(X)$.  Also, let $\mathbb
Q(X)(\alpha_1,\ldots,\alpha_t)(y)$ be a finite extension of
$\mathbb Q(X)(\alpha_1,\ldots,\alpha_t)$ defined by the equation
$q_i(X,y) = 0$.  It follows that $y$ is algebraic over $\mathbb
Q(X)$, and upon clearing denominators, let $\tilde q_i(X,Y) \in
\mathbb Q[X,Y]$ be such that $\tilde q_i(X,y) = 0$.

Since both $q_i(X,Y)$ and $\tilde q_i(X,Y)$ have $y$ as a root and
since $q_i(X,Y)$ is irreducible, it follows that $q_i(X,Y)$
divides $\tilde q_i(X,Y)$.  As both $\tilde q_i(X,Y)$ and
$q_i(X,Y)$ are in $\mathbb Q^*[X,Y]$, Gauss's Lemma \cite[p.
181]{Lang} gives us that, \[h_i(X,Y) = \tilde q_i(X,Y)/q_i(X,Y)
\in \mathbb Q^*[X,Y].\]  Let \[\tilde g(X,Y) = \prod\limits_{i =
1}^r {\tilde q_i(X,Y)} \in \mathbb Q[X,Y]. \] We claim that
$\tilde g(X,Y)$ is our desired polynomial.  But indeed, for all $p
\in A \cap \mathbb Q^*$, \[\tilde g(p,\phi_j(p)) = \prod\limits_{i
= 1}^r {q_i(p,\phi_j(p))h_i(p,\phi_j(p))} = 0.\]

Finally, let $A$ and $R$ be as in the statement of Theorem
\ref{moregenthm}, and let $W$ be the graph of $\phi_j: A
\rightarrow R$, which is a semialgebraic set (over $\mathbb Q$).
Consider the sentence,
\[\forall p \ \forall a \ \neg((p,a) \in W) \
\vee \ \tilde g(p,a) = 0,\] which is valid in any structure
satisfying the axioms of RCF (again, since $\mathbb Q$ embeds in
any real closed field). By the above argument, it is a true
statement for $\mathbb Q^*$, and therefore, by completeness, it is
also true for $R$. This completes the proof. $\QED$

As a final remark, we note that the fussiness in the proof of
Theorem \ref{moregenthm} was necessary to avoid division by zero
after applying the substitution homomorphism with $X \mapsto p$
and $Y \mapsto \phi_j (p)$.




\begin{thebibliography}{5}
\bibitem{Al}
D. J. Aldous (1990).
A random walk construction of uniform spanning trees and uniform labeled trees.
{\em SIAM Journal on Discrete Mathematics}, {\bf 3(4)}, 450--465.


\bibitem{Bo}
A. Broder (1989). Generating random spanning trees.
In {\em 30th Annual Symposium on Foundations of Computer Science}, 442--447.


\bibitem{CLO}
D. Cox, J. Little and D. O'Shea, Donal (1997).
{\em Ideals, varieties, and algorithms.
An introduction to computational algebraic geometry and commutative algebra.}
Second edition. Undergraduate Texts in Mathematics. Springer-Verlag, New York, 1997.

\bibitem{CS}
N. Chomsky and M. P. Sch\"{u}tzenberger (1963).
The algebraic theory of context-free languages.
In {\em Computer programming and formal systems}, 118--161, North-Holland, Amsterdam.

\bibitem{El}
P. Elias (1972). The efficient construction of unbiased random sequence, {\em Ann. Math. Stat.}
{\bf 43}, 865--870.

\bibitem{HU}
J. E. Hopcroft and J. D. Ullman (1979).
{\em Introduction to automata theory, languages, and computation.}
Addison-Wesley Series in Computer Science. Addison-Wesley Publishing Co., Reading, Mass.

\bibitem{Gr}
A. Grzegorczyk (1955). Computable functionals. {\em Fund. Math.}
 {\bf 42}, 168--202.

\bibitem{GH}
P.W. Glynn and S. Henderson.
Nonexistence of a class of variate generation schemes. 
{\em Operations Research Letters}, to appear. 

\bibitem{HLP}
G.H. Hardy, J.E. Littlewood, and G. P\'{o}lya (1959).
{\em Inequalities},
Cambridge University Press, Cambridge.


\bibitem{KO}
M. S. Keane , G. L. O'Brien (1994).
A Bernoulli factory, {\em ACM Transactions on Modeling and Computer Simulation} {\bf 4}, Issue 2.

\bibitem{Po}
G. P\'{o}lya (1928) \"{U}ber positive Darstellung von Polynomen
Vierteljschr, 
{\em Naturforsch. Gez. \"{Z}urich}
{\bf 73}, 141--145. In {\em Collected papers} {\bf 2} (1974) MIT press, 309--313.

\bibitem{Pe}
Y. Peres (1992). Iterating von Neumann's Procedure for Extracting Random Bits,
{\em Ann. Stat.}, {\bf 20}, 590--597.

\bibitem{PR}
M. B. Pour-El and J. I. Richards (1988). {\em Computability in Analysis and Physics,} Springer-Verlag.

\bibitem{PR2}
V. Powers and B. Reznick (2002).
A new bound for Polya's Theorem with applications to polynomials positive on polyhedra,
to appear in {\em MEGA 2000 proceedings, J. Pure Applied Algebra}.

\bibitem{PW1}
J. G. Propp and D. B. Wilson (1996).
Exact sampling with coupled Markov chains and applications to statistical mechanics.
{\em Random Structures and Algorithms}, {\bf 9(1,2)}, 223--252.

\bibitem{PW2}
J. G. Propp and D. B. Wilson (1998).
How to get a perfectly random sample from a generic Markov chain and generate a random spanning tree of a directed graph.
{\em Journal of Algorithms} {\bf 27}, 170--217.


\bibitem{Tur}
A. Turing (1936). On computable numbers, with an application to the Entscheidungsproblem.
{\em Proc. London Math. Soc.}, {\sc Ser.\ 2} {\bf 42}, 230--265.


\bibitem{vN}
J. von Neumann (1951).
Various techniques used in connection with random digits.
{\em Applied Math Series} {\bf 12}, 36--38.

\end{thebibliography}

\begin{thebibliography}{9}

\bibitem{Boch} J. Bochnak, M. Coste, and M. F. Roy, \emph{Geometrie Algebrique Reelle},
Springer Verlag, 1986.

\bibitem{Dickman} M. Dickman, \emph{Applications of Model Theory to Real Algebraic
geometry, in Methods in Mathematical Logic}, Springer Verlag,
1985.

\bibitem{Haskell} D. Haskell, A. Pillay, and C. Steinhorn, \emph{Model Theory, Algebra, and Geometry},
Cambridge University Press, 2000.

\bibitem{Lang} S. Lang, \emph{Algebra -3rd ed}, Addison-Wesley Publishing Company, New
York, 1993.

\bibitem{Marker} D. Marker.  \emph{Model Theory: an Introduction}, Springer Verlag, 2002.

\bibitem{Marker2} D. Marker. \emph{Introduction to the Model Theory of Fields},
http://www.math.uic.edu/\textasciitilde marker/.

\end{thebibliography}
\end{document}